\documentclass{amsart}
\usepackage{graphicx}
\usepackage{latexsym}
\usepackage{hyperref}
\usepackage{amsfonts}
\usepackage[all]{xy}
\usepackage{amssymb, mathrsfs, amsfonts, amsmath}
\usepackage{amsbsy}
\usepackage{amsfonts}
\newtheorem{theorem}{Theorem}

\numberwithin{equation}{section}

\begin{document}

\title {\bf A modified Yamabe invariant and a Hopf conjecture }
\author {Ezio de Araujo Costa}

\date{}

\maketitle

\begin{center}{\bf Abstract}

\end{center}
In this paper we define the bi-orthogonal sectional curvature and we present two modified Yamabe invariants for compact 4-dimensional manifolds. In particular we obtained a relationship between one of these invariants and a Hopf conjecture.
\\

\noindent {\bf Mathematic subject classifications (2000): 53C25, 53C24. \noindent
}

\section{\bf Introduction and Statments of results}

\begin{center}{\bf The bi-orthogonal sectional curvature}
\end{center}
  Let $M = M^n$  be a compact manifold of dimension $n\geq 4$ and denotes by ${\bf Met}(M)$ the set of Riemannian metrics on $M$. For each metric $g \in {\bf Met}(M)$, let $s$ be the scalar curvature of $M$ in metric $g$ and denotes by $K$ the sectional curvature of this metric. For each $x\in M$, let $P_1, P_2$ two mutually orthogonal two dimensional subspace of tangent space $T_xM$. We call of {\it bi-orthogonal sectional curvature} ($K^\perp$) relative to $P_1$ and $P_2$ (in $x\in M$) the number given by $$K^\perp(P_1,P_2) = \frac {K(P_1) + K(P_2)}{2} \eqno [1.1]$$
If $n=4$, we write $$K^\perp(P) = \frac {K(P) + K(P^\perp)}{2}. \eqno [1.2]$$.

This notion of curvature was used in [7] and [6] by W. Seaman and M. H. Noronha, respectively.
In particular, Seaman proved an extension of part of the Sphere Theorem for $n\geq 5$ (see Theorem 0.3 in [7]).
Now, let  $M$  be a 4-dimensional Riemannian manifold and consider the following functions on $M$:

$$k_1^\perp = \textsl{min} \ \{K^\perp(P), \ P \ \textmd{a 2-plan in } T_xM \}, \eqno[1.3]$$

$$k_3^\perp = \textsl{max} \ \{K^\perp(P), \  P \ \textmd{a 2-plan in } T_xM \}, \eqno[1.4]$$

$$k_2^\perp = \frac{s}{4} - k_1^\perp - k_3^\perp. \eqno [1.5]$$

 In four dimension, we have the following features of the bi-orthogonal sectional curvature:
\\
\\
I) $M$ is an Einstein 4-manifold if and only if $K^\perp(P) = K(P)$, for each $x\in M$ and for all 2-plain $P \subset T_xM$.
\\
\\
II) A 4-manifold $M$ is conformally flat if and only if $k_i^\perp(x) = \frac{s(x)}{12}$ ($i = 1, 2, 3$), for all $x\in M$, where $s$ is the scalar curvature of $M$.
\newpage
{\bf Remark }- Property I is a well known characterization of Einstein 4-manifold. Property II follows from a  criterion of Kulkarni.
\\
\\
\begin{center} {\bf  An extension of the Sphere Theorem in four dimension}
\end{center}
If the bi-orthogonal of a manifold $M$ satisfies $1/4 \leq K^\perp \leq 1$ then $M$ has nonnegative isotropic curvature (see [7]). Using  the classification of compact manifolds with non negative isotropic curvature curvature (see [1] and [8]) we obtained an extension of a theorem of Seaman (see Theorem 0.3 in [ 7]):
\\
\\
\begin{theorem}
\label{thm1}
  Let $M$ be a  compact oriented 4-manifold such that $1/4 \leq K^\perp \leq 1$ and let the  $\widetilde{M}$ be the universal covering of $M$.
\begin{enumerate}
 \item If $M$ is irreducible then $\widetilde{M}$ is diffeomorphic to sphere $\mathbb{S}^4$, $M$ is biholomorphic  to  complex projective space $\mathbb{CP}^2$ or $M$ is diffeomorphic to a connected sum $\mathbb{S}^4 \sharp m\sharp \mathbb{RP}^4 \sharp (\mathbb{R} \times \mathbb{S}^3)/G_1 \sharp ... \sharp (\mathbb{R} \times \mathbb{S}^3)/G_n$, where $m = 0$, $m =1$ and the $G_i$ are discrete subgroup of the isometry group of $\mathbb{R}\times \mathbb{S}^3$.
 \item
  If $M$ is reducible then $\widetilde{M}$ is isometric to a Riemannian product $\mathbb{R} \times N^2$, where $N^2$ is diffeomorphic to sphere $\mathbb{S}^2$.
\end{enumerate}
\end{theorem}
 \begin{center}
 {\bf  The bi-orthogonal sectional curvature and the Weyl tensor in four dimension}
\end{center}
The Weyl tensor $W$ of an oriented 4-manifold $M$ has the following decomposition: $W = W^+\bigoplus W^-$, where $W^\pm$ are the self-dual and anti-self-dual parts of the tensor $W$, respectively.
Let $\Lambda ^2$ be the space of two-forms $\varphi$ in $M$ and let $* :\Lambda ^2 \rightarrow \Lambda^2$ be the star operator of Hodge. Then $ \Lambda^+ \bigoplus \Lambda^-$, where $\Lambda^\pm = \{\varphi \in \Lambda^2; \ *\varphi = \pm\varphi \}$. The Weyl tensor has the decomposition  $W = W^+ \bigoplus W^-$, where $W^\pm : \Lambda^\pm \rightarrow \Lambda^\pm $ are self-adjoint operators with free traces. $W^\pm$ are called the self-dual and anti-self-dual parts of the Weyl tensor $W$ of $M$, respectively. The matrix of the curvature operator $\Re$ of $M$ takes the form
\[
 \left(\begin{array}{clcr}
W^+ + \frac{s}{12}I_{\Lambda^-} & \ \ B \\
B^* & \ W^- + \frac{s}{12}I_{\Lambda^-}
\end{array}\right),
\]
where  $B:\Lambda^- \rightarrow \Lambda^+$, $\mid B \mid^2 = \mid Ric - \frac{s}{4}\mid^2$ and $Ric$ is the Ricci operator of $M$.
\\
\\
Denote by $w_1^\pm \leq w_2^\pm \leq w_3^\pm$ the eigenvalues of $W^\pm$, respectively.
\\
\\
Let $x\in M$ and consider $X, Y$ orthonormal in tangent space $T_xM$. Then a simple and unitary two-form $\varphi = X\wedge Y$ can be uniquely written as $\varphi = \varphi^+ + \varphi^-$, where $\varphi^\pm \in \Lambda^\pm$ and $\mid \varphi^\pm \mid^2 = \frac{1}{2}$. The sectional curvature $K(\varphi)$ is given by
$$K(\phi) = \frac{s}{12} + \langle \varphi^+, W^+(\varphi^+) \rangle + \langle \varphi^-, W^-(\varphi^-) \rangle + 2\langle \varphi^+, B\varphi^-\rangle \eqno[1.6]$$

If $\varphi^-$ is replaced by $ -\varphi^-$, we have

$$K(\varphi^\perp) = \frac{s}{12} + \langle \varphi^+, W^+(\varphi^+) \rangle + \langle \varphi^-, W^-(\varphi^-) \rangle  - 2\langle \varphi^+, B\varphi^-\rangle, \eqno[1.7]$$
Where $\varphi^\perp = \varphi^+  - \varphi^-$.
Adding [1.6]  to [1.7] :
$$\frac{K(\varphi^\perp) + K(\varphi)}{2} = \frac{s}{12} + \langle \varphi^+, W^+(\varphi^+) \rangle + \langle \varphi^-, W^-(\varphi^-) \rangle$$
Using [1.3],

$$k_{1}^\perp = \ \frac{s}{12} \ + \textsl{min} \ \{ \langle \varphi^+, W^+(\varphi^+) \rangle + \langle \varphi^-, W^-(\varphi^-) \rangle, \mid \varphi^\pm \mid^2 = 1/2\}  =$$
$$\frac{s}{12} \ + \textsl{min} \ \{ \langle \varphi^+, W^+(\varphi^+) \rangle, \ \mid \varphi^+ \mid^2 =1/2\} + \textsl{min} \ \{ \langle \varphi^-, W^-(\varphi^-) \rangle, \ \mid \varphi^- \mid^2 =1/2\}.$$
\\
\\
By Propositions 2.1 in [6] there  exists an orthonormal basis $\{ P_1, P_2, P_3, P_4\}$ of $\Lambda^2$ such that each $P_i$ is the form $X_i \wedge Y_i $, where $X_i, Y_i \in T_xM$. Moreover, in accord with the Proposition 2.5 in [6] is easy see that the eigenvectors of $W^\pm$ are
$\frac {\sqrt{2}}{2}( P_i \pm P_i^\perp)$, respectively. Then

$$k_{1}^\perp = \frac{s}{12} + \frac{w_1^+ + w_1^-}{2}, \eqno [1.8]$$
where $w_1^\pm$ are the smallest eigenvalues of $W^\pm(x)$, respectively.
\\
\\
Similarly and in view of [1.4] we have

$$k_{3}^\perp = \frac{s}{12} + \frac{w_3^+ + w_3^-}{2}, \eqno [1.9]$$
where $w_3^\pm$ are the largest eigenvalues of $W^\pm(x)$, respectively.
\\
\\
Since that $w_2^\pm = -w_1^\pm - w_3^\pm$, respectively, we can uses [1.5] and obtain
$$k_{2}^\perp =  \frac{s}{12} + \frac{w_3^+ + w_3^-}{2}. \eqno [1.10]$$

An extension of the Yamabe Problem  was considered by M. Itoh in [4] and more recently, M. Listimg (see chapter 2, section 2.2 in [5]) and B-L Chen and X-P Zhu (see [2]) obtained important results on this topic.
In our article we use some results of [2] for two modified Yamabe invariants.
\begin{center}{\bf A modified Yamabe invariant and the isotropic curvature. }
\end{center}
 Let $M$ be a compact oriented 4-manifold and let ${\bf Met}(M)$ be the set of Riemannian metrics on $M$. If $g \in {\bf Met}(M)$,
  let $[g] = \{\widetilde{g} = u^2g ; u\in C^\infty(M), u > 0\}$ and let $dV_{\widetilde{g}}$ be the volume element in metric $\widetilde{g}$.
  \\
  \\
  The Yamabe constant of the metric $g$ is given by
$$Y (M, g):= \textsf{inf} \ \{ \frac{1}{\sqrt{V_{\widetilde{g}}}} \int_M \widetilde{s} dV_{\widetilde{g}}, \ \widetilde{g} \in [g] \},$$ where $\widetilde{s}$ is the scalar curvature of $M$ in metric $\widetilde{g}$.
\\
\\
The Yamabe invariant of $M$ is given by
$$Y (M):= sup \  \{ Y_g(M) \}.$$
Now let

$$Y^\perp (M, g):= \textsf{inf} \ \{ \frac{1}{\sqrt{V_{\widetilde{g}}}} \int_M [24k_1^\perp - \widetilde{s}]dV_{\widetilde{g}}, \ \widetilde{g} \in [g] \}, \eqno [1.11]$$
where $ \widetilde{k}_1^\perp$ is the bi-orthogonal sectional curvature of $M$ given by [1.3] in metric $\widetilde{g}$. Consider the following modified Yamabe invariant:
$$Y^\perp (M) := sup \ \{ Y^\perp(M, g),  g\in {\bf Met}(M)\}. \eqno [1.12]$$
Our next result is the following
\newpage

\begin{theorem}
\label{thm2}
  Let $M$ be a compact oriented 4-manifold with Riemannian metric $g$. If $Y^\perp (M, g) \geq 0$
   then we have
\begin{enumerate}

  \item $M$ is diffeomorphic to a connected sum $\mathbb{S}^4\sharp m\sharp \mathbb{RP}^4 \sharp (\mathbb{R} \times \mathbb{S}^3)/G_1 \sharp ... \sharp (\mathbb{R} \times \mathbb{S}^3)/G_n$, where $m = 0$ or $1$, $i \geq 0$ and the $G_i$ are discrete subgroup of the isometry group of $\mathbb{R}\times \mathbb{S}^3$ or
 \item $(M, g)$ is  conformal to a complex projective space $\mathbb{CP}^2$ with the Fubini-Study metric or a finite cover is conformal to a Riemannian product  $\mathbb{S}_{c_1}^2 \times \mathbb{T}^2$, where $\mathbb{S}_{c} ^2$ is a sphere with constant sectional curvature $c$ and $\mathbb{T}^2$ is a flat torus.
 \end{enumerate}
 \end{theorem}

 The conditions of the Theorem 2 imply that $M$ admits a metric with non negative isotropic curvature.
 \\
 \\
 \begin{center}{\bf A modified Yamabe invariant and a Hopf conjecture }
 \end{center}
 Let $M$ be a compact oriented 4-manifold and $g\in {\bf Met}(M)$
and consider
$$Y_1^\perp (M, g):= \textsf{inf} \ \{ \frac{1}{\sqrt{V_{\widetilde{g}}}} \int_M k_1^\perp dV_{\widetilde{g}}, \ \widetilde{g} \in [g] \}.  \eqno [1.13]$$ We have another  modified Yamabe invariant:
$$Y_1^\perp (M) := sup \ \{ Y_1^\perp(M, g), g \in {\bf Met}(M)\}.\eqno [1.14]$$

 Recall the Hopf conjecture:
  \\
  \\
  {\bf $\mathbb{S}^2 \times \mathbb{S}^2$ no admits a Riemannian metric with positive sectional curvature.}
   \\
   \\
  Then we can formulate the following question:
   \\
   \\
   {\bf  $\mathbb{S}^2 \times \mathbb{S}^2$  admits a metric with positive bi-orthogonal sectional curvature ?}
   \\
   \\
   With respect to this question we have:
\\
\\
\begin{theorem}
\label{thm3}
 Let $M$ be a compact oriented 4-manifold. Then we have
\begin{enumerate}
 \item $M$ has a metric $g$  with  $k_1^\perp > 0$ if and only if $Y_{1}^\perp (M) > 0$.
 \item $Y_1^\perp (M) \leq Y(M) \leq Y(\mathbb{S}^4)$
, where $Y(M)$ is the Yamabe invariant of $M$. In particular, if $Y_1^\perp (M) = Y(\mathbb{S}^4)$ then $M$ is conformal to the standard sphere $\mathbb{S}^4$.
\item If $M$ has a metric $g$  with  $k_1^\perp \geq 0$ and scalar curvature $s > 0$ then $$8\pi^2\chi < max \{ \int_{M}s^2dV_g + 16\pi^2, \frac{5}{24} \int_{M}s^2dV_g \}, $$
where $\chi$ is the Euler characteristic of $M$.
\end{enumerate}
\end{theorem}

 {\bf Corollary 4}
\\
 (1) {\it  Let $g_{can}$ the canonical metric of product of spheres $\mathbb{S}^2 \times \mathbb{S}^2$ and let $g \in [g_{can}]$. Then $g$ no has $k_1^\perp > 0$ on $\mathbb{S}^2 \times \mathbb{S}^2$.}
\\
\\
 (2) {\it If $Y_1^\perp (\mathbb{S}^2 \times \mathbb{S}^2) \leq  0$ then the Hopf conjecture is true.}
  \\
  \\
  (3) {\it Let $g$ a Riemannian metric on $\mathbb {S}^2 \times \mathbb{S}^2$ with scalar curvature $s$. If $\int s^2dV_g \leq \frac{ 768\pi^2}{5}$ then  $g$ no has $k_1^\perp \geq 0$.}

\section{\bf Proof of results}
\begin{center}{\bf Proof of Theorem 1}
\end{center}
 Let $M$ be a  compact oriented 4-manifold such that $1/4 \leq K^\perp \leq 1$. In accord with Seaman [7], $M$ has non negative isotropic curvature. Consider $M$ is irreducible. By main result of Seshadri [8] we have only the following possibilities:
 \\
 \\
 i) $M$ admits a metric with positive isotropic curvature. In this case, follow from the main  theorem of [1] that
$M$ is diffeomorphic to a connected sum $\mathbb{S}^4 \sharp m\sharp \mathbb{RP}^4 \sharp (\mathbb{R} \times \mathbb{S}^3)/G_1 \sharp ... \sharp (\mathbb{R} \times \mathbb{S}^3)/G_n$, where $m = 0$, $m =1$ and the $G_i$ are discrete subgroup of the isometry group of $\mathbb{R}\times \mathbb{S}^3$.
 \\
 \\
 ii) $M$ is isometric to a irreducible locally symmetric space. Then $M$ is an Einstein space with positive isotropic curvature and so $M$ is isometric to a sphere $\mathbb{S}^4$.
\\
\\
 iii) $M$ is biholomorphic to a complex projective space $\mathbb{CP}^2$.
 \\
\\
Now, assumes that $M$ is reducible. Since that $M$ has positive scalar curvature, then the universal covering $\widetilde{M}$ is  isometric to $\mathbb{R} \times M_2^3$ or isometric to $M_1^2 \times M_2^2$. But $M_1^2 \times M_2^3$ has $k_1^\perp = 0$. So, $\widetilde{M}$ is isometric to a Riemannian product $\mathbb{R} \times N^2$, where $M_2^3$ has positive sectional curvature and in this case $M_2^3$ is diffeomorphic to sphere $\mathbb{S}^3$.
This finish the proof of Theorem 1.
\\
\\
For  proof of Theorems 2 and 3  we need a lemma. For this let $M$ be a compact oriented 4-manifold $M$ with metric $g$ and scalar curvature $s$ and consider the functions $f(W) = 2s - 24k_1^\perp = -12(w_1^+ + w_1^-) \geq 0$ (see [1.8] )and $f_1(W) = s - 12k_1^\perp  = -6(w_1^+ + w_1^-) \geq 0$  , where $w_1^\pm$ are the smallest eigenvalues of $W^\pm$, respectively.  In notation of [2 , eq. (2.2) and (2.5)] we have  $Y_{f}(M, [g]) = Y^\perp(M,g)$ and $Y_{f_1} (M, [g]) = Y_1^\perp (M, g)$, where $Y^\perp(M,g)$ and $Y_1^\perp (M, g)$ are given by [1.11] and [1.13], respectively.
\\
\\
 In accord with the lemma 2.1 and 2.2 in [2] we have:
 \\
 \\
{\bf Lemma 5}- {\it Let $M$ be a compact oriented 4-manifold with metric $g$.
\\
\\
 i) There exists $ \widetilde{g}, \overline{g} \in [g]$ such that $Y^\perp (M, g) =  24\widetilde{k}_1^\perp - \widetilde{s} = $ constant and $Y_1^\perp (M, g) = 12\overline{k}_1^\perp  = $ constant, respectively where
$ \widetilde{s}$ is the scalar curvature in metric $\widetilde{g}$ and $\widetilde{k}_1^\perp$ and $\overline{k}_1^\perp$ are the smallest bi-orthogonal sectional curvatures in metrics $\widetilde{g}$ and $\overline{g}$, respectively.
\\
\\
 ii) If $Y^\perp (M, g) > 0$ or $Y_1^\perp (M, g) > 0$ then there exists $ \widetilde{g} \in [g]$ such that $  24\widetilde{k}_1^\perp - \widetilde{s} > 0 $ or $\widetilde{k}_1^\perp > 0$, respectively.}
\\
\\
\begin{center}{\bf Proof of Theorem 2}
\end{center}
Let $M$ be a compact oriented 4-manifold with Riemannian metric $g$ and $Y^\perp (M, g) \geq 0$. Initially assumes that $Y^\perp (M, g) > 0$. By Lemma 5)ii there exists $ \widetilde{g} \in [g]$ such that $ 24\widetilde{k}_1^\perp - \widetilde{s} > 0 $. Then (see [1.8]) we have $$ w_1^\pm \geq w_1^+ + w_1^- >  -\frac{s}{12}$$

Since that $w_3^\pm  \leq -2w_1^\pm$, respectively then $\frac{s}{6} - w_3^\pm > 0$ an this proves that $(M, \widetilde{g})$ has positive isotropic curvature. By main result of [1], $M$ is diffeomorphic to a connected sum $\mathbb{S}^4\sharp m\sharp \mathbb{RP}^4 \sharp (\mathbb{R} \times \mathbb{S}^3)/G_1 \sharp ... \sharp (\mathbb{R} \times \mathbb{S}^3)/G_n$, where $m = 0$ or $1$, $i \geq 0$ and the $G_i$ are discrete subgroup of the isometry group of $\mathbb{R}\times \mathbb{S}^3$
 \\
Now, consider  $Y^\perp (M, g) = 0$. By Lemma 5)i there exists
  $ \widetilde{g}$ such that $Y^\perp (M, g) = 24\widetilde{k}_1^\perp - \widetilde{s} = 0$.
   In this case $(M, \widetilde{g})$ has nonnegative isotropic curvature. Assumes that $M$ is irreducible. By main result of Seshadri [8] we have only the following possibilities:
   \\
 \\
 i) $M$ admits a metric with positive isotropic curvature. In this case, follow from the main  theorem of  that
$M$ is diffeomorphic to a connected sum $\mathbb{S}^4 \sharp m\sharp \mathbb{RP}^4 \sharp (\mathbb{R} \times \mathbb{S}^3)/G_1 \sharp ... \sharp (\mathbb{R} \times \mathbb{S}^3)/G_n$, where $m = 0$, $m =1$ and the $G_i$ are discrete subgroup of the isometry group of $\mathbb{R}\times \mathbb{S}^3$.
 \\
 \\
 ii) $M$ is isometric to a irreducible locally symmetric space. Then $M$ is an Einstein space with positive isotropic curvature and so $M$ is isometric to a sphere $\mathbb{S}^4$.
\\
\\
 iii) $(M, g)$ is conformal to a complex projective space $\mathbb{CP}^2$ with the Fubini-Study metric.
 \\
 \\
Now, assumes that $M$ is reducible.  Then $(M, g)$ is either conformal to a finite cover Riemannian product $\mathbb{S}_{c_1}^2 \times \mathbb{S}_{c_2}^2$, $\mathbb{S}_{c}^3 \times \mathbb{R}$ or $\mathbb{S}_{c_1}^2 \times \mathbb{T}^2$, where $\mathbb{S}_{c}$ is a sphere with constant sectional curvature $c$ and $\mathbb{T}^2$ is a flat torus. But $\mathbb{S}_{c_1}^2 \times \mathbb{S}_{c_2}^2$ has $k_1 = 0 \neq s/24$. This finish the proof of Theorem 2.

\begin{center}{\bf Proof of Theorem 3}
\end{center}
Theorem 3)(1) and 3)(2) are consequences of  Corollary 2.2.3 and  equation (2.2) in [5] , respectively.
\\
  Let $M$ with metric $g$,  $k_1^\perp \geq 0$ and scalar curvature $s > 0$. In accord Theorem 1.1 in [9], if $\chi$ is the Euler characteristic of $M$ then  $M$ is  isometric to a sphere or $$8\pi^2(\chi - 2) <  \int_{M}\mid W\mid ^2dV_g, \eqno$$
 where $W$ is the Weyl tensor of $(M, g)$. Note that we have the inequalities :
 $$\mid W \mid ^2 = \mid W^+ \mid^ 2 + \mid W^-\mid^ 2 \leq 6\mid w_1^+\mid^2 + 6\mid w_1^-\mid^2 \leq 6(w_1^+ + w_1^-)^2$$
 Using [1.8], we obtain  $\mid W \mid^2 \leq 6(\frac{s}{6} - 2k_1)^2 \leq \frac{s^2}{6}$ and so $8\pi^2\chi < \frac{1}{6} \int_{M} s^2 dV_g.$
In any case have that $8\pi^2\chi < \frac{1}{6} \int_{M} s^2 dV_g + 16\pi^2$
 On the other side, if $(M, g)$ is a compact oriented 4-manifold with metric $g$ and scalar curvature $s$ then the Euler characteristic of $M$ satisfy  $$8\pi^2\chi = \int_{M} [\mid W\mid^2 + \frac{s^2}{24} - \frac{1}{2}\mid B\mid^2]dV_g,$$ where $B = 0$ if and only if $g$ is an Einstein metric.
 Using previous inequalities we have $8\pi^2\chi \leq \frac{5}{24} \int_{M} s^2 dV_g.$ In particular, if $8\pi^2\chi = \frac{5}{24} \int_{M} s^2 dV_g$ then we can see that $B = 0$ and $g$ is an Einstein metric with nonnegative sectional curvature  and positive scalar curvature. By Lemma 2 in [3], $8\pi^2\chi < \frac{5}{24}\int_M s^2dV_g$.
 In any case have that $8\pi^2\chi < \frac{5}{24} \int_{M} s^2 dV_g.$ This finish the proof the Theorem 3
\begin{center}{\bf Proof of Corollary 4 }
\end{center}
  (1) Let $g_{can}$ be the canonical metric of product of spheres $\mathbb{S}^2 \times \mathbb{S}^2$ and let $g \in [g_{can}]$. Consider $\widetilde{k}_1^\perp $ and $k_1^\perp $ the smallest bi-orthogonal curvatures of the metrics
  $g$ and
  $g_{can}$, respectively. Let $g = u^2g_{can}$. Then $12\widetilde{k}_1^\perp = u^{-3}[-6\Delta_{g}u + 12k_1^\perp]$ (see eq. (2.3) in [2]). Since that $k_1^\perp = 0$ we have $\int_{M}\widetilde{k}_1^\perp = 0$ and this proves that $g$ no has $\widetilde{k}_1^\perp > 0. $
  \\
  \\
  (2) Corollary 4)(2) is consequence of Theorem 3)(1).
  \\
  \\
  (3) Let $g$ be a Riemannian metric on $\mathbb {S}^2 \times \mathbb{S}^2$ with $k_1^\perp \geq 0$. Assumes that  the scalar curvature $s$ of $g$ satisfy $\int_{M} s^2dV_g \leq \frac{ 768\pi^2}{5}$. Then by Corollary 3)(3)3, $\chi < 4$ which contradicts the fact of that $\chi(\mathbb {S}^2 \times \mathbb{S}^2) = 4$

Author's address:

Mathematics Department,  Federal University of Bahia,

zipcode: 40170110- Salvador -Bahia-Brazil

Author's email : ezio@ufba.br
\end{document}